\documentclass[12pt]{article}
\usepackage{graphicx}
\usepackage{amsmath,amsthm,amssymb,enumerate}
\usepackage{euscript,mathrsfs}
\usepackage[left=2cm,right=2cm,top=3.5cm,bottom=3.5cm]{geometry}
\usepackage{color}
\catcode`\@=11 \@addtoreset{equation}{section}

\catcode`\@=12

\allowdisplaybreaks

\newtheorem{Theorem}{Theorem}[section]
\newtheorem{Proposition}[Theorem]{Proposition}
\newtheorem{Lemma}[Theorem]{Lemma}
\newtheorem{Corollary}[Theorem]{Corollary}

\theoremstyle{definition}
\newtheorem{Definition}[Theorem]{Definition}

\newtheorem{Remark}[Theorem]{Remark}

\newcommand{\bTheorem}[1]{
\begin{Theorem} \label{T#1} }
\newcommand{\eT}{\end{Theorem}}

\newcommand{\bProposition}[1]{
\begin{Proposition} \label{P#1}}
\newcommand{\eP}{\end{Proposition}}

\newcommand{\bLemma}[1]{
\begin{Lemma} \label{L#1} }
\newcommand{\eL}{\end{Lemma}}

\newcommand{\bCorollary}[1]{
\begin{Corollary} \label{C#1} }
\newcommand{\eC}{\end{Corollary}}

\newcommand{\bRemark}[1]{
\begin{Remark} \label{R#1} }
\newcommand{\eR}{\end{Remark}}

\newcommand{\bDefinition}[1]{
\begin{Definition} \label{D#1} }
\newcommand{\eD}{\end{Definition}}

\newcommand{\bfphi}{\boldsymbol{\varphi}}

\newcommand{\bFormula}[1]{
\begin{equation} \label{#1}}
\newcommand{\eF}{\end{equation}}

\newcommand{\Ov}[1]{\overline{#1}}

\newcommand{\DC}{C^\infty_c}

\newcommand{\vr}{\varrho}
\newcommand{\vre}{\vr_\ep}

\newcommand{\vue}{\vu_\ep}

\newcommand{\vu}{\vc{u}}

\newcommand{\vc}[1]{{\bf #1}}
\newcommand{\Ome}{\Omega_\ep}

\newcommand{\Div}{{\rm div}_x}
\newcommand{\Grad}{\nabla_x}

\newcommand{\dx}{\,{\rm d} {x}}

\newcommand{\dt}{\,{\rm d} t }

\newcommand{\intO}[1]{\int_{\Omega} #1 \ \dx}

\newcommand{\intOe}[1]{\int_{\Omega_\ep} #1 \ \dx}

\newcommand{\intRN}[1]{\int_{R^3} #1 \ \dx}

\newcommand{\ep}{\varepsilon}

\newcommand\be{\begin{equation}}
\newcommand\ee{\end{equation}}
\newcommand\ba{\begin{equation}\begin{aligned}}
\newcommand\ea{\end{aligned}\end{equation}}

\definecolor{Cgrey}{rgb}{0.85,0.85,0.85}
\definecolor{Cblue}{rgb}{0.50,0.85,0.85}
\definecolor{Cred}{rgb}{1,0,0}
\definecolor{fancy}{rgb}{0.10,0.85,0.10}

\newcommand\Cbox[2]{%
    \newbox\contentbox%
    \newbox\bkgdbox%
    \setbox\contentbox\hbox to \hsize{%
        \vtop{
            \kern\columnsep
            \hbox to \hsize{%
                \kern\columnsep%
                \advance\hsize by -2\columnsep%
                \setlength{\textwidth}{\hsize}%
                \vbox{
                    \parskip=\baselineskip
                    \parindent=0bp
                    #2
                }%
                \kern\columnsep%
            }%
            \kern\columnsep%
        }%
    }%
    \setbox\bkgdbox\vbox{
        \color{#1}
        \hrule width  \wd\contentbox %
               height \ht\contentbox %
               depth  \dp\contentbox
        \color{black}
    }%
    \wd\bkgdbox=0bp%
    \vbox{\hbox to \hsize{\box\bkgdbox\box\contentbox}}%
    \vskip\baselineskip%
}

\date{}

\begin{document}


\title{Homogenization of a non--homogeneous heat conducting fluid}

\author{Eduard Feireisl
\thanks{The research of E.F.~leading to these results has received funding from the
European Research Council under the European Union's Seventh
Framework Programme (FP7/2007-2013)/ ERC Grant Agreement
320078. The Institute of Mathematics of the Academy of Sciences of
the Czech Republic is supported by RVO:67985840.} \and Yong Lu \and Yongzhong Sun \thanks{The research of Y. Sun was partially supported by NSF of China under Grant 11571167.}
}

\date{\today}

\maketitle

\bigskip

\centerline{Institute of Mathematics of the Academy of Sciences of the Czech Republic}

\centerline{\v Zitn\' a 25, CZ-115 67 Praha 1, Czech Republic}

\centerline{Email: feireisl@math.cas.cz}

\centerline{and}

\centerline{Department of Mathematics, Nanjing University}

\centerline{22 Hankou Road, Gulou District, 210093 Nanjing, China}
\centerline{ Email: luyong@nju.edu.cn}

\centerline{and}

\centerline{Department of Mathematics, Nanjing University}

\centerline{22 Hankou Road, Gulou District, 210093 Nanjing, China}
\centerline{ Email: sunyz@nju.edu.cn}

\bigskip

\begin{abstract}

We consider a non--homogeneous incompressible and heat conducting fluid confined to a 3D domain perforated by tiny holes.
The ratio of the diameter of the holes and their mutual distance is critical, the former being equal to $\ep^3$, the latter proportional to $\ep$, where $\ep$ is a small parameter. We identify the asymptotic limit
for $\ep \to 0$, in which the momentum equation contains a friction term of Brinkman type determined uniquely by the viscosity and geometric properties of the perforation. Besides the inhomogeneity of the fluid, we allow the viscosity and the heat conductivity 
coefficient to depend on the temperature, where the latter is determined via the Fourier law with homogenized (oscillatory) heat conductivity coefficient that is different for the fluid and the solid holes. To the best of our knowledge, this is the first result 
in the critical case for the inhomogenous heat--conducting fluid.

\end{abstract}



\section{Introduction}
\label{i}

We study the motion of a non--homogeneous, incompressible viscous and heat conducting fluid contained in a bounded spatial domain perforated by a system
of tiny holes. The mass density $\vr = \vr(t,x)$, the velocity $\vu = \vu(t,x)$ and the temperature $\Theta = \Theta(t,x)$ satisfy
a variant of the Navier--Stokes--Boussinesq
system proposed by Chandrasekhar \cite{CHAN} (see also Ligni\`eres \cite{LIGN}) :
\begin{equation} \label{i1}
\partial_t \vr + \Div (\vr \vu) = 0, \ \Div \vu = 0,
\end{equation}
\begin{equation} \label{i2}
\begin{split}
\partial_t (\vr \vu) + \Div (\vr \vu \otimes \vu) + \Grad P &= \Div \mathbb{S}(\Theta, \Grad \vu) - \Theta \Grad F,\\
\mathbb{S}(\Theta,  \Grad \vu) &= \mu(\Theta) \left( \Grad \vu + \Grad^{\rm T} \vu \right),
\end{split}
\end{equation}
\begin{equation} \label{i3}
- \Div \left( \kappa \Grad \Theta \right) = \Grad F \cdot \vu.
\end{equation}
Here the last equation can be seen as a quasi--static (high P\' eclet number) approximation of the conventional heat equation
\[
\partial_t (\vr \Theta) + \Div (\vr \Theta \vu) - \Div \left( \kappa \Grad \Theta \right) = \Grad F \cdot \vu,
\]
where $F$ is the gravitational potential.
We refer to \cite[Chapter 4, Section 4.3]{FENO6} for a rigorous derivation of system (\ref{i1}--\ref{i3}) in the spatially homogeneous case $\vr \equiv 1$.

The fluid is contained in a bounded domain $\Omega_\ep \subset R^3$, on the boundary of which the velocity obeys the no-slip condition
\begin{equation} \label{i4}
\vu|_{\partial \Omega_\ep} = 0.
\end{equation}
Extending $\vu$ to be zero outside $\Omega_\ep$ we may therefore assume that the equation of continuity (\ref{i1}) is satisfied in the whole space $R^3$.
Similarly, we suppose $\Omega_\ep \subset \Omega$ whereas (\ref{i3}) is satisfied in $\Omega$, with
\[
\kappa = \kappa_\ep (x) = \left\{ \begin{array}{l} \kappa_f > 0 \ \mbox{if}\ x \in \Omega_\ep, \\ \\
\kappa_s > 0 \ \mbox{if}\ x \in \Omega \setminus \Ov{\Omega}_\ep \end{array} \right.
\]
where, in general, we allow $\kappa_f \ne \kappa_s$. For definiteness, we prescribe the homogeneous Dirichlet boundary conditions for the temperature,
\begin{equation} \label{i5}
\Theta |_{\partial \Omega} = 0.
\end{equation}

\subsection{Perforated domain}
\label{iS1}

We now introduce the perforated domain under consideration. Let $0<\ep<1$ be a small parameter. We suppose
\[
\Omega_\ep = \Omega \setminus \bigcup_{k = 1}^{K(\ep)} T_{k, \ep},
\]
where
\begin{equation} \label{i6}
T_{k,\ep} = x_k + \ep^3 \Ov{U}_{k,\ep},\ k=1,2,\dots, K(\ep),
\end{equation}
with ${\rm dist} [x_i, x_j] > c\, \ep \ \mbox{whenever}\ i \ne j, \ {\rm dist} [x_i, \partial \Omega] > c \,\ep,
\ \mbox{for some constant}\ c > 0 \ \mbox{independent of}\ \ep$. By a normalization process, we may assume $c=1$.

Here $\{ U_{k,\ep} \}_{\ep > 0, k = 1, \dots, K(\ep)}$ are assumed to be uniformly $C^{2 + \nu}$ simply connected domains satisfying
\be\label{size-holes}
\left\{ | x | < \frac{1}{2} \right\} \subset U_{k,\ep} \subset \Ov{U}_{k,\ep} \subset \left\{ | x | < \frac{3}{4} \right\}
\ \mbox{for any}\ \ep, \ k.
\ee

Thus possible spatial configuration of the holes $T_{k,\ep}$ includes the so--called critical case, the holes being of radius $\ep^3$ with their
mutual distance proportional to $\ep$, cf. Allaire \cite{Allai4} among others. The assumptions \eqref{i6}--\eqref{size-holes} imposed on the distribution of holes guarantee the holes are pairwise disjoint. Note that no periodicity of the holes is a priori assumed.

Finally, for the sake of simplicity, we suppose that $\partial \Omega$ is smooth, of class $C^{2 + \nu}$.  We use $C$ to denote a universal constant whose value is independent of $\ep$.

\subsection{Weak solutions}

We consider weak solutions to problem (\ref{i1}--\ref{i5}) emanating from the initial data
\begin{equation} \label{i7}
\vr(0, \cdot) = \vr_{0, \ep}, \ \vu(0, \cdot) = \vu_{0, \ep} ,
\end{equation}
and belonging to the regularity class
\begin{equation} \label{i8}
\begin{split}
\vr &\in C([0,T]; L^1(\Omega)), \ 0 < \underline{\vr} \leq \vr \leq \Ov{\vr} \ \mbox{a.a. in}\ (0,T) \times \Omega,\\
\Theta &\in L^\infty(0,T; W^{1,2}_0(\Omega)), \ \Theta(t, \cdot) \in C^\nu (\Ov{\Omega}),\
\| \Theta(t, \cdot) \|_{C^\nu (\Ov{\Omega})} \leq C \ \mbox{for a.a.}\ t \in (0,T),\ \nu > 0,\\
\vu &\in L^\infty(0,T; L^2(\Omega; R^3)) \cap L^2(0,T; W^{1,2}_0(\Omega; R^3)).
\end{split}
\end{equation}

The equations (\ref{i1}), (\ref{i2}), (\ref{i3}) will be interpreted in the weak sense. More specifically,
\begin{equation} \label{i9}
\int_0^T \intRN{ \left[ \vr \partial_t \varphi + \vr \vu \cdot \Grad \varphi \right] } \dt = - \intRN{ \vr_{0,\ep} \varphi(0, \cdot)}
\end{equation}
for any $\varphi \in C^1_c([0,T) \times R^3)$, where $\vu \equiv 0$ outside $\Omega_\ep$;
\begin{equation} \label{i10}
\begin{split}
\Div \vu &= 0 \ \mbox{a.a. in}\ (0,T) \times \Omega;\\
&\int_0^T \intOe{ \left[ \vr \vu \cdot \partial_t \bfphi + \vr \vu \otimes \vu : \Grad \bfphi \right] } \dt \\
&= \int_0^T \intOe{ \mathbb{S}(\Theta, \Grad \vu) : \Grad \bfphi } + \int_0^T \intO{ \Theta \Grad F \cdot \bfphi } \dt
- \intOe{ \vr_{0,\ep} \vu_{0,\ep} \cdot \bfphi(0, \cdot) }
\end{split}
\end{equation}
for any $\bfphi \in C^1_c([0,T) \times \Ome; R^3))$, $\Div \bfphi = 0$;
\begin{equation} \label{i11}
\intO{ \kappa_\ep \Grad \Theta (\tau, \cdot) \cdot \Grad \phi } = \intO{ \vu (\tau, \cdot) \cdot \Grad F \phi }
\end{equation}
for a.a. $\tau \in (0,T)$ and $\phi \in C^1_c (\Omega)$.

In addition, we suppose that the energy inequality
\begin{equation} \label{i12}
\begin{aligned}
\intOe{ \frac{1}{2} \vr |\vu|^2 (\tau, \cdot) }
+ \int_0^\tau \intO{ \kappa_\ep |\Grad \Theta|^2 } \dt + \int_0^\tau \intOe{ \frac{\mu(\Theta)}{2} |\Grad \vu + \Grad^t \vu |^2 } \dt\\
\leq \intOe{ \frac{1}{2} \vr_{0,\ep} |\vu_{0,\ep} |^2 }
\end{aligned}
\end{equation}
holds for a.a. $\tau \in (0,T)$.

In view of the DiPerna--Lions theory \cite{DL}, and the anticipated regularity of $\vr$, $\vu$ stated in (\ref{i8}), the weak formulation (\ref{i9}) implies its renormalized variant
\begin{equation} \label{i13}
\int_0^T \intRN{ \left[ b(\vr) \partial_t \varphi + b(\vr) \vu \cdot \Grad \varphi \right] } \dt = - \intRN{ b(\vr_{0,\ep}) \varphi(0, \cdot)}
\end{equation}
for any $\varphi \in C^1_c([0,T) \times R^3)$ and any $b \in  C((0,\infty))$ (actually any Borel $b$), due to the lower and upper bound restriction of $\vr$.

\medskip

We remark that,  for any fixed $\ep>0$, the existence of renormalized weak solutions can be derived following the nowadays well understood argument in Lions's book \cite{Lions-Incom}.

\subsection{Main result}

Let $K \subset \left\{ |x| < 1 \right\}$ be a compact set. We define the matrix
\[
\mathbb{C}_{i,j} (K) = \int_{\{ |x| < 1 \} \setminus K} \Grad \vc{v}^i : \Grad \vc{v}^l \ \dx, \ i,j=1, \dots 3,
\]
where $\vc{v}^i$ is the unique solution of the \emph{model problem}
\begin{equation} \label{i14}
- \Delta \vc{v}^i + \Grad q^i = 0, \ \Div \vc{v}_i = 0 \ \mbox{in}\  \left\{ |x| < 1 \right\} \setminus K,\
\vc{v}^i|_{\partial K} = \vc{e}^i, \ \vc{v}^i|_{\{ |x| = 1 \}} = 0.
\end{equation}
Here $\{ \vc{e}^i \}_{i=1}^3$ denotes the standard orthogonal basis of the vector space $R^3$.
Note that $\vc{v}^i$ is the unique minimizer of the Dirichlet integral
\[
\intRN{ |\Grad \vc{v} |^2 } \ \mbox{over the set}\ \{\vc{v} \in W^{1,2}(R), \ \Div \vc{v} = 0, \ \vc{v}|_K = \vc{e}^i, \ \vc{v}|_{\{ |x| \geq 1 \}} = 0\}.
\]

Finally, we suppose that there is a positive definite symmetric matrix field $\mathbb{D} \in L^\infty(\Omega; R^{3 \times 3}_{{\rm sym}})$ such that
\begin{equation} \label{i15}
\lim_{\ep \to 0} \sum_{T_{k,\ep} \subset B} \mathbb{C}_{i,j} (\ep^3 \Ov{U}_{k,\ep}) = \int_B \mathbb{D}_{i,j}(x) \ \dx
\ \mbox{for any Borel set}\ B \subset \Omega,
\end{equation}
where $U_{k,\ep}$ are related to $T_{k,\ep}$ via (\ref{i6}). Note that the limit (\ref{i15}) exists in the spatially periodic case with holes of
uniform shape studied in the nowadays classical papers by Allaire \cite{Allai4}, \cite{Allai3}. Other relevant examples can be found in Desvillettes, Golse, and Ricci \cite{DesGolRic},
or Marchenko, Khruslov \cite{MarKhr}.

We are ready to formulate our main result:

\begin{Theorem} \label{T1}
Let $\{ \Ome \}_{\ep > 0}$ be a family of perforated domains specified in Section \ref{iS1}, where the asymptotic distribution of holes
satisfies (\ref{i5}). Let the initial data be given such that
\begin{equation} \label{i16}
\begin{split}
\vr_{0,\ep} &\in L^\infty(R^3), \ \vr_{0,\ep}(x) = \vr_s > 0 - \mbox{a positive constant} \ \mbox{for}\ x \in R^3 \setminus \Ov{\Omega}_\ep\\
0 &< \underline{\vr} \leq \vr_{0,\ep}(x) \leq \Ov{\vr}, \ \vr_{0,\ep} \to \vr_0 \ \mbox{in}\ L^1(\Omega);
\end{split}
\end{equation}
\begin{equation} \label{i17}
\Div \vu_{0,\ep} = 0, \ \vu_{0,\ep} = 0 \ \mbox{in}\ R^3 \setminus \Ov{\Omega}_\ep, \ \vu_{0,\ep} \to \vu_0 \ \mbox{in}\ L^2(\Omega; R^3);
\end{equation}
Finally, suppose that $\Grad F \in L^\infty(\Omega; R^3)$ and that $\mu = \mu(\Theta)$ is a positive continuous function of $\Theta$.

Let $( \vre, \vue, \Theta_\ep)$ be a weak solution of the problem (\ref{i1}--\ref{i3}), (\ref{i4}), (\ref{i5}), (\ref{i7}). Then, up to a subsequence, we have
\[
\begin{split}
\vre &\to \vr \in C([0,T]; L^1(\Omega)),\ 0 < \underline{\vr} \leq \vr_{\ep}(t, x) \leq \Ov{\vr}, \\
\Theta_\ep &\to \Theta \ \mbox{in}\ L^q((0,T) \times \Omega)\ \mbox{for any}\ 1 \leq q < \infty,
\ \mbox{and weakly-(*) in}\ L^\infty(0,T; W^{1,2}_0(\Omega)), \\
\vue &\to \vu \ \mbox{in}\ L^2((0,T) \times \Omega; R^3) \ \mbox{and weakly in}\ L^2(0,T; W^{1,2}_0(\Omega; R^3)),
\end{split}
\]
where $(\vr, \vu, \Theta)$ is a weak solution of the problem
\[
\partial_t \vr + \Div (\vr \vu) = 0, \ \Div \vu = 0,
\]
\[
\begin{split}
\partial_t (\vr \vu) + \Div (\vr \vu \otimes \vu) + \mu(\Theta) \mathbb{D} \vu + \Grad P &= \Div \mathbb{S}(\Theta, \Grad \vu) + \Theta \Grad F,\\
\mathbb{S}(\Theta,  \Grad \vu) &= \mu(\Theta) \left( \Grad \vu + \Grad \vu^t \right),
\end{split}
\]
\[
- \Div \left( \kappa_f \Grad \Theta \right) = \Grad F \cdot \vu,
\]
in $(0,T) \times \Omega$, satisfying the boundary conditions (\ref{i4}), (\ref{i5}), and the initial conditions $(\vr_0, \vu_0)$.

\end{Theorem}

The rest of the paper is devoted to the proof of Theorem \ref{T1}. It is worth noting that the limit process includes in fact \emph{two} homogenization procedures: the first one in the momentum equation due to the domain perforation, the second one in the heat equation due to the spatial oscillations
of the heat conductivity coefficient. The two processes interact via the temperature dependent viscosity coefficient $\mu$. Besided the nowadays standard
homogenization technique developed in the pioneering paper by Allaire \cite{Allai4}, our method leans essentially on the uniform estimates of the H\" older norm
of the temperature $\Theta_\ep$. To the best of our knowledge, this is the first result concerning the critical case for the inhomogeneous 
heat conducting fluid. It is worth noting that the Brikman type term in the asymptotic limit is independent of the density of the fluid, cf. 
the nowadays classical paper of Cioranescu and Murat \cite{CioMur2}, \cite{CioMur1} concerning the background of this extra term.

The paper is organized as follows. In Section \ref{p}, we derive some preliminary estimates that follow directly from the renormalized formulation and the available energy bounds, in particular we derive the uniform bounds on the H\" older norm of the temperature. The homogenization process in the momentum
equation is performed in Section \ref{H}. Finally, the limit passage is completed in Section \ref{a}. To conclude, let us remark that, in contrast with the bulk of the available homogenization literature almost exclusively focused on stationary problems, the evolutionary setting requires essential modifications of
the limit process.

\section{Preliminaries - uniform estimates}
\label{p}

We start with uniform bounds for $\vr$ and $\Theta$. Using hypothesis (\ref{i16}) we can take
\[
b(\vr) = [\vr - \Ov{\vr}]^+, \ b(\vr) = - [\vr - \underline{\vr}]^-
\]
as test functions in the renormalized equation (\ref{i13}) to deduce
\begin{equation} \label{p1}
0 < \underline{\vr} \leq \vre(t,x) \leq \Ov{\vr} \ \mbox{for a.a.}\ (t,x)
\end{equation}
uniformly in $\ep \to 0$. Next, using the lower bound for (\ref{p1}) for $\vr$, we
we deduce from the energy inequality (\ref{i12}), combined with (\ref{i17}),
that
\begin{equation} \label{p4}
{\rm ess}\sup_{t \in [0,T]} \left\| \vue (t, \cdot) \right\|_{L^2(\Omega; R^3)} + \int_0^T \| \Grad \Theta_\ep \|^2_{L^2(\Omega)} \leq C
\end{equation}
Note that $\vue \equiv 0$ outside $\Omega_\ep$.

Now, seeing that $\Theta_\ep$ solves for a.a. fixed time the elliptic equation (\ref{i11}), with the diffusion coefficient
\[
0 < \min \{ \kappa_s, \kappa_f \} \leq \kappa \leq \max \{ \kappa_s, \kappa_f \},
\]
we may use the standard elliptic theory, see e.g.  Ladyzhenskaya, Uralceva \cite[Chapter 3, Theorem 12.1]{LADUR}, to obtain the estimate
\begin{equation} \label{p44}
{\rm ess}\sup_{t \in (0,T)} \| \Theta_\ep (t, \cdot) \|_{C^\nu(\Ov{\Omega})} \leq C
\end{equation}
for a certain $\nu > 0$. It is important that the bound in (\ref{p44}) depends solely on $\kappa_s$, $\kappa_f$, and the constant in (\ref{p4}), specifically on the norm of the initial data.
In particular,
\begin{equation} \label{p2}
 - \Ov{\Theta} \leq \Theta_\ep (t,x) \leq \Ov{\Theta} \ \mbox{for a.a.}\ (t,x).
\end{equation}
Going back to the energy balance (\ref{i12}) and using the positivity of $\mu$ on the range $[-\Ov{\Theta}, \Ov{\Theta}]$
we may infer that
\begin{equation} \label{p4E}
\int_0^T \| \Grad \vue \|^2_{L^2(\Omega; R^3)} \leq C.
\end{equation}

Now, using (\ref{p1}), (\ref{p4E}), the renormalized equation (\ref{i13}) and hypothesis (\ref{i16}), we get
\begin{equation} \label{p5}
\begin{split}
\vre &\to \vr \ \mbox{in}\ C_{\rm weak}([0,T]; L^q(\Omega)) \ \mbox{for any}\ 1 < q < \infty,\\
\vue &\to \vu \ \mbox{weakly-(*) in}\ L^\infty(0,T; L^2(\Omega; R^3))
\ \mbox{and weakly in}\ L^2(0,T; W^{1,2}_0(\Omega;R^3)),
\end{split}
\end{equation}
passing to suitable subsequences as the case may be. In addition, the standard Aubin--Lions argument yields immediately that $\vr$, $\vu$ satisfy
(\ref{i9}). Finally, by DiPerna--Lions theory \cite{DL}, the same equation holds in the renormalized sense (\ref{i13}). In particular, it can be shown that
\[
\| \vre \|^2_{L^2(\Omega)} \to \| \vr \|^2_{L^2(\Omega)} \ \mbox{in}\ C[0,T];
\]
whence
\begin{equation} \label{p6}
\vre \to \vr \ \mbox{in}\ C([0,T]; L^q(\Omega)) \ \mbox{for any}\ 1 \leq q < \infty.
\end{equation}

\section{Homogenization}
\label{H}

We start with the elliptic problem associated to the momentum equation (\ref{i2}):
\begin{equation} \label{H1}
- \Div \left[ \mu(\Theta) \left( \Grad \vc{U}_\ep + \Grad^t \vc{U}_\ep \right) \right] + \Grad P_\ep = \vc{f}_\ep,\
\Div \vc{U}_\ep = 0 \ \mbox{in} \ \Omega_\ep ,\ \vc{U}_\ep |_{\partial \Omega_\ep} = 0.
\end{equation}

For a given $\Theta \in C^\nu (\Ov{\Omega})$ and $\vc{f}_\ep \in W^{-1,2}(\Ome; R^3)$, problem (\ref{H1}) admits
a weak solution $\vc{U}_\ep$, $P_\ep$, unique in the class
\be\label{Ue-Pe-class}
\vc{U}_\ep \in W^{1,2}_0(\Ome; R^3), \ P_\ep \in L^2(\Ome),\ \intOe{ P_\ep } = 0,
\ee
such that the equations in \eqref{H1} are satisfied in the weak sense: for any $\phi \in C_c^\infty(\Omega_\ep)$ and any $\boldsymbol{\varphi} \in C_c^\infty(\Omega_\ep;R^3)$, there holds
\begin{equation} \label{H1-weak1}
\int_{\Omega_\ep} \vc{U}_\ep \cdot \Grad \phi \,\dx =0,
\end{equation}
and
\begin{equation} \label{H1-weak2}
\int_{\Omega_\ep} \mu(\Theta) \left( \Grad \vc{U}_\ep + \Grad^t \vc{U}_\ep \right): \Grad \boldsymbol{\varphi}  - P_\ep \Div \boldsymbol{\varphi}\,\dx = \langle \vc{f}_\ep,  \boldsymbol{\varphi}\rangle_{W^{1,2}_0(\Omega_\ep)}.
\end{equation}

We remark that the solution can be obtained as the minimizer of the functional
\[
\vc{U} \mapsto \intRN{ \frac{1}{2} \mu(\Theta) | \Grad \vc{U} + \Grad^t \vc{U} |^2 }  -  \left< \vc{f}_\ep \cdot \vc{U} \right>
\]
over the space of functions
\[
\left\{ \vc{U} \in W^{1,2}(R^3; R^3) \ \Big| \ \Div \vc{U} = 0, \ \vc{U}|_{R^3 \setminus {\Omega}_\ep  =  0} \right\}.
\]

Our goal in this section is to show the following result.

\begin{Proposition} \label{P1}

Let $\{ \Omega_\ep \}_{\ep > 0}$ be a family of domains satisfying the same hypotheses as Theorem \ref{T1}.
Suppose that
\be\label{P1-ass}
\Theta \in C^\nu (\Ov{\Omega}), \ \limsup_{\ep \to 0 } \left\| \vc{f}_\ep  - \vc{f} \right\|_{W^{-1,2}(\Ome; R^3)} \leq M
\ee
for some $\vc{f}\in W^{-1,2}(\Omega; R^3)$ independent of $\ep$ and some $M \geq 0$.

Let $\vc{U}_\ep$, $P_\ep$ be the unique (weak) solution of problem (\ref{H1}). Then, up to the zero extension and a substraction of subsequence, there holds
\ba\label{Ue-Pe-weak-limit}
\vc{U}_\ep \to \vc{U} \ \mbox{weakly in} \ W^{1,2}_0(\Omega; R^3),\  {P}_\ep \to P
\ \mbox{weakly in}\ L^2(\Omega) \ \mbox{as}\ \ep \to 0,
\ea
where $\vc{U}$, $P$ is the solution of the problem
\ba\label{U-P-equaions}
- \Div \left[ \mu (\Theta) \left( \Grad \vc{U} + \Grad^t \vc{U} \right) \right] + \mu(\Theta)
\mathbb{D} \vc{U} + \Grad P = \vc{f} + \vc{r},\
\Div \vc{U} = 0 \ \mbox{in} \ \Omega, \ \vc{U}|_{\partial \Omega} = 0
\ea
for some $\vc{r} \in W^{-1,2}_0(\Omega; R^3)$ satisfying
\[
\| \vc{r} \|_{W^{-1,2}_0 (\Omega; R^3)} \leq C M.
\]

\end{Proposition}

The rest of this section is devoted to proving Proposition \ref{P1}. This is done in the following subsections step by step by employing similar arguments as in \cite{FNN16}, where the main idea goes back to \cite{DesGolRic}.

\medskip

We recall the the following pointwise and integral estimates of the solution $(\vc{v}^i,q^i)$ to the model problem \eqref{i14}. The proof follows from the proof of Lemma 4.1 in \cite{FNN16}.
\begin{Lemma}\label{model-pt}

Let $(\vc{v}^i,q^i)$ be a solution to \eqref{i14} with $K \subset B(0,r) \subset B(0,d) \subset B(0,1)$. Then there holes the estimates
\be\label{model-ptwise}
\left| \partial^\alpha \vc{v}^i \right| \leq C \frac{r}{|x|^{1+|\alpha|}}, \quad \left| q^i   \right| \leq C \frac{r}{|x|^{2}}, \quad \forall x\in B(0,1) \setminus B(0,r),
\ee
where $\alpha \in \mathbb{N}^3, \ |\alpha| \leq 2,$ and

\ba\label{model-L2}
&\int_{B(0,d)} |\vc{v}^i |^2\,\dx \leq C r^2 d,  \quad  \int_{B(0,d)} | \Grad \vc{v}^i |^2\,\dx \leq C r,   \\
&\int_{B(0,d)} |q^i |^2\,\dx \leq C r,
\ea

\end{Lemma}

\subsection{Uniform estimates}

Since $\Theta \in C^\nu (\Ov{\Omega})$, it admit a lower and upper bound. By the assumption that $\mu(\Theta)$ is positive and continuous function in $\Theta$, we have that for some positive constants $\underline\mu$ and $\bar\mu$,
\be\label{mu(theta)-low-high}
0< \underline\mu \leq \mu(\Theta) \leq \bar\mu <\infty.
\ee

By a density argument, the weak formulation \eqref{H1-weak2} is satisfied for any $\boldsymbol{\varphi} \in W^{1,2}_0(\Omega_\ep;R^3)$. By \eqref{Ue-Pe-class}, we can take the solution $\vc{U}_\ep$ itself to be a test function in \eqref{H1-weak2}  and obtain
\ba\label{Ue-est1}
\int_{\Omega_\ep} \mu(\Theta) \left| \Grad \vc{U}_\ep + \Grad^t \vc{U}_\ep \right|^2 \,\dx  &= 2  \langle \vc{f}_\ep, \vc{U}_\ep \rangle_{W^{1,2}_0(\Omega_\ep)} \\
&\leq 2 \|\vc{f}_\ep\|_{W^{-1,2}(\Omega_\ep)} \|\vc{U}_\ep\|_{W_0^{1,2}(\Omega_\ep)} \leq 2 (\|\vc{f}\|+M) \|\vc{U}_\ep\|_{W_0^{1,2}(\Omega_\ep)} .
\ea

By \eqref{mu(theta)-low-high}-\eqref{Ue-est1}, applying Korn's inequality and Poincar\'e's inequality gives
\ba\label{Ue-est2}
\|\vc{U}_\ep\|_{W_0^{1,2}(\Omega_\ep)}^2 & \leq C \|\Grad \vc{U}_\ep\|_{L^{2}(\Omega_\ep)}^2 \leq C \|\Grad \vc{U}_\ep + \Grad^t \vc{U}_\ep\|_{L^{2}(\Omega_\ep)}^2 \\
& \leq C \int_{\Omega_\ep} \mu(\Theta) \left| \Grad \vc{U}_\ep + \Grad^t \vc{U}_\ep \right|^2 \,\dx   \leq C (\|\vc{f}\|_{W^{-1,2}(\Omega_\ep)} + M) \|\vc{U}_\ep\|_{W_0^{1,2}(\Omega_\ep)}.
\ea
This implies the uniform estimate:
\ba\label{Ue-est-f}
\sup_{0<\ep<1} \|\vc{U}_\ep\|_{W_0^{1,2}(\Omega_\ep)} \leq C  (\|\vc{f}\|_{W^{-1,2}(\Omega_\ep)} + M).
\ea

Theorem 2.3 in \cite{DFL17} applies to the setting of perforated domains in this paper. As a result, there exits a linear uniform bounded Bogovskii type operator
\[
\mathcal{B}_\ep : L_0^{2}(\Omega_\ep) \to W_0^{1,2}(\Omega_\ep; R^3),
\]
such that for any $f\in L_0^{2}(\Omega_\ep)$,
\ba\label{pt-bog}
\Div \mathcal{B}_\ep(f) =f \ \mbox{in} \ \Omega_\ep,\quad \|\mathcal{B}_\ep(f)\|_{W_0^{1,2}(\Omega_\ep; R^3)}\leq C\, \|f\|_{L^2(\Omega_\ep)},
\ea
for some constant $C$ independent of $\ep$.

Since $P_\ep \in L_0^2(\Omega_\ep)$ which is the collection of $L^2(\Omega_\ep)$ functions with zero mean value, we have $\mathcal{B}_\ep(P_\ep) \in W_0^{1,2}(\Omega_\ep;R^3)$. Taking $\mathcal{B}_\ep(P_\ep)$ as a test function in the weak formulation \eqref{H1-weak2} implies
\ba\label{Pe-est1}
 \int_{\Omega_\ep}   |P_\ep|^2 \,\dx = \int_{\Omega_\ep} \mu(\Theta) \left( \Grad \vc{U}_\ep + \Grad^t \vc{U}_\ep \right): \Grad \mathcal{B}_\ep(P_\ep) \,\dx -  \langle \vc{f}_\ep,  \mathcal{B}_\ep(P_\ep)\rangle_{W^{1,2}_0(\Omega_\ep)}.
\ea
Together with \eqref{Ue-est-f} and \eqref{pt-bog}, we obtain from \eqref{Pe-est1} that
\ba\label{Pe-est2}
\| P_\ep \|_{L^2(\Omega_\ep)}^2 \leq  C  (\|\vc{f}\|_{W^{-1,2}(\Omega_\ep)} + M) \|\mathcal{B}_\ep(P_\ep)\|_{W^{1,2}_0(\Omega_\ep)} \leq C (\|\vc{f}\|_{W^{-1,2}(\Omega_\ep)} + M) \| P_\ep \|_{L^2(\Omega_\ep)},
\ea
which implies
\ba\label{Pe-est-f}
\sup_{0<\ep<1} \|P_\ep\|_{L^{2}(\Omega_\ep)} \leq C  (\|\vc{f}\|_{W^{-1,2}(\Omega_\ep)} + M).
\ea

Hence, by the uniform estimates in \eqref{Ue-est-f} and \eqref{Pe-est-f}, up to the zero extensions and a substraction of subsequence, there holds
\be\label{Ue-Pe-limit}
\vc{U}_\ep \to \vc{U} \ \mbox{weakly in} \ W^{1,2}_0(\Omega; R^3),\  {P}_\ep \to P
\ \mbox{weakly in}\ L^2(\Omega) \ \mbox{as}\ \ep \to 0.
\ee

By the divergence free property of $\vc{U}_\ep$, we have $\Div \vc{U} =0$. It is left to prove that the limit $(\vc{U},P)$ solves the Brinkman type equations in \eqref{U-P-equaions}. This is done in the next subsection.

\subsection{Decompositions}

Let $\chi$ be a function satisfying
\be\label{chi-def}
\chi\in C_c^\infty(-1,1), \quad \chi=1 \ \mbox{on}  \ \big[-\frac{3}{4},\frac{3}{4}\big],  \quad 0\leq \chi \leq 1, \ 0\leq \chi' \leq 4.
\ee
Define the cut-off function $\phi_{k,\ep}$ near each hole $T_{k,\ep}$ by
\be\label{phi-k-e-def}
\phi_{k,\ep}(x) : = \chi \left( \frac{|x - x_k|}{\ep^3} \right).
\ee

By the assumptions \eqref{i6}--\eqref{size-holes} of the distribution of holes, there holds
$$
\phi_{k,\ep}(x) = 1 \ \mbox{on} \ T_{k,\ep}, \quad \phi_{k,\ep} \phi_{k',\ep} =0 \ \mbox{whence $k\neq k'$}.
$$

Let $\vc{v}^i_{\ep}$ be the solution to the model problem \eqref{i14} with $K = \overline U_{k,\ep}$. For any $\boldsymbol\varphi = (\boldsymbol\varphi^i)_{i=1,2,3} \in C_c^\infty(\Omega;R^3)$, we define $\boldsymbol\varphi^{\ep}\in W_0^{1,2}(\Omega_\ep;R^3)$ as
\be\label{Restriction-direct}
\boldsymbol\varphi^{\ep} := \boldsymbol\varphi - \boldsymbol\varphi^{\ep}_1 - \boldsymbol\varphi^{\ep}_2,
\ee
with
\ba\label{varphi-e-1-2}
\boldsymbol\varphi^{\ep}_1(x) & :=  \sum_{k=1}^{K(\ep)} \left(\boldsymbol\varphi(x) - \boldsymbol\varphi(x_k) \right) \phi_{k,\ep}(x),\\
\boldsymbol\varphi^{\ep}_2(x) & :=  \sum_{k=1}^{K(\ep)} \sum_{i=1}^3  \boldsymbol\varphi^i(x_k)  \vc{v}^i_{\ep} ( x - x_k ) \phi_{k,\ep}(x).
\ea
Give the above definition, it is immediately to find
\be
\boldsymbol\varphi^{\ep}  = 0, \quad \mbox{on} \ \bigcup_{k=1}^{ K(\ep)} T_{k,\ep}.
\ee
So there does holds $\boldsymbol\varphi^{\ep}\in W_0^{1,2}(\Omega_\ep;R^3)$ with
$$
\|\boldsymbol\varphi^{\ep}\|_{W_0^{1,2}(\Omega_\ep;R^3)} \leq C \|\boldsymbol\varphi\|_{W_0^{1,2}(\Omega;R^3)}.
$$
Moreover, similarly as Lemma 5.2 in \cite{FNN16}, a direct calculation gives
\ba\label{varphi-e-1-2-limit}
&\boldsymbol\varphi^{\ep}_1 \to 0, \ \mbox{strongly in} \  W^{1,2}_0(\Omega;R^3), \ \mbox{as $\ep\to 0$},\\
&\boldsymbol\varphi^{\ep}_2 \to 0, \ \mbox{weakly in} \  W^{1,2}_0(\Omega;R^3), \ \mbox{as $\ep\to 0$}.
\ea

Let $(\vc{U},P)$ be the limit we obtained in \eqref{Ue-Pe-limit}. For any given $\kappa>0$, there exits $\vc{U}^0 \in C_c^\infty(\Omega;R^3)$ such that
\be\label{U-Ukappa}
\vc{U}  =  \vc{U}^0 + \vc{U}^\kappa, \quad \|\vc{U}^\kappa\|_{W^{1,2}(\Omega)} \leq \kappa.
\ee
As in \eqref{Restriction-direct} and \eqref{varphi-e-1-2}, we consider the decomposition
\be\label{U-U0-12}
 \vc{U}^0 =  \vc{U}^0_\ep + \vc{U}_{\ep,1}^0 +  \vc{U}_{\ep,2}^0,
\ee
where $\vc{U}_{\ep,1}^0$ and $\vc{U}_{\ep,2}^0$ are defined in the same manner as in \eqref{varphi-e-1-2} and satisfy the same convergence results as in \eqref{varphi-e-1-2-limit}.

We thus consider the decomposition of $\vc{U}_\ep$ as
\be\label{Ue-dec}
\vc{U}_\ep = \vc{U}^0_\ep + \vc{U}_\ep^r,
\ee
 where $\vc{U}^0_\ep$ comes from the decomposition \eqref{U-U0-12}. It is crucial to study the property of the remainder $\vc{U}_\ep^r = \vc{U}_\ep - \vc{U}^0_\ep \in W_0^{1,2}(\Omega_\ep;R^3).$ Due to the fact
$$
\vc{U}_\ep \to \vc{U} =  \vc{U}^0 + \vc{U}^\kappa, \ \vc{U}^0_\ep \to \vc{U}^0, \ \mbox{weakly in} \ W^{1,2}_0(\Omega;R^3), \quad \mbox{as $\ep\to 0$},
$$
it is straightforward to obtain that
\be\label{Ur-limit0}
\vc{U}_\ep^r \to  \vc{U} - \vc{U}^0 =  \vc{U}^\kappa, \ \mbox{weakly in} \ W^{1,2}_0(\Omega;R^3), \quad \mbox{as $\ep\to 0$.}
\ee
Then the Rellich-Kondrachov compact embedding theorem implies, up to a substraction of subsequence, that
\be\label{Ur-limit1}
\vc{U}_\ep^r \to    \vc{U}^\kappa, \ \mbox{strongly in} \ L^{q}(\Omega;R^3), \quad \forall q \in [2, 6),\quad \mbox{as $\ep\to 0$}.
\ee
A consequence is that
\be\label{Ur-limit2}
\limsup_{\ep \to 0} \|\vc{U}_\ep^r\|_{L^q(\Omega)} \leq C(q) \|\vc{U}^\kappa\|_{L^q(\Omega)} \leq C(q) \|\vc{U}^\kappa\|_{W^{1,2}_0(\Omega)} \leq C(q)\kappa, \quad  \forall q \in [2, 6).
\ee
Moreover, it can be shown that
\be\label{Ur-limit3}
\limsup_{\ep \to 0} \|\nabla \vc{U}_\ep^r\|_{L^2(\Omega)}  \leq h(\kappa) \to 0,\  \mbox{as} \ \kappa \to 0.
\ee
In order to estimate $\|\nabla \vc{U}_\ep^r\|_{L^2(\Omega)}$, choosing $\vc{U}_\ep^r$ as a test function in the weak formulation  \eqref{H1-weak2} implies
 \ba \label{H1-weak2-Ur}
& \int_{\Omega_\ep} \mu(\Theta) \left( \Grad \vc{U}_\ep^r + \Grad^t \vc{U}_\ep^r \right): \Grad \vc{U}_\ep^r =  - \int_{\Omega_\ep} \mu(\Theta) \left( \Grad \vc{U}_\ep^0 + \Grad^t \vc{U}_\ep^0 \right): \Grad \vc{U}_\ep^r \,\dx\\
& \qquad + \int_{\Omega_\ep}   P_\ep \Div \vc{U}_\ep^r\,\dx + \langle \vc{f}_\ep, \vc{U}_\ep^r  \rangle_{W^{1,2}_0(\Omega_\ep)}\\
&\quad  =  - \int_{\Omega_\ep} \mu(\Theta) \left( \Grad \vc{U}^0 + \Grad^t \vc{U}^0 \right): \Grad \vc{U}_\ep^r \,\dx  + \int_{\Omega_\ep}   P_\ep \Div \vc{U}_\ep^r\,\dx + \langle \vc{f}_\ep, \vc{U}_\ep^r  \rangle_{W^{1,2}_0(\Omega_\ep)}\\
& \qquad+ \int_{\Omega_\ep} \mu(\Theta) \left( \Grad \vc{U}_{\ep,1}^0 + \Grad^t \vc{U}_{\ep,1}^0 \right): \Grad \vc{U}_\ep^r \,\dx + \int_{\Omega_\ep} \mu(\Theta) \left( \Grad \vc{U}_{\ep,2}^0 + \Grad^t \vc{U}_{\ep,2}^0 \right): \Grad \vc{U}_\ep^r \,\dx.
\ea
Then starting from \eqref{H1-weak2-Ur}, by \eqref{Ur-limit0}--\eqref{Ur-limit2}, together with the strong convergence of  $\vc{U}_{\ep,1}^0$  and weak convergence of $\vc{U}_{\ep,2}^0$, by using the property of $\vc{v}^i_{\ep}$ as the solution to the model problem \eqref{i14} with $K = \overline U_{k,\ep}$ (see Lemma \ref{model-pt}), a similar argument as the proof of Lemma 5.1 in \cite{FNN16} implies our desired result in \eqref{Ur-limit3}.

\subsection{Limit equations}

Now we deduce the equations satisfied by the limit couple $(\vc{U},P)$. Let $\boldsymbol\varphi = (\boldsymbol\varphi^i)_{i=1,2,3} \in C_c^\infty(\Omega;R^3)$ and let  $\boldsymbol\varphi^{\ep}\in W_0^{1,2}(\Omega_\ep;R^3)$ be defined as in \eqref{Restriction-direct}. Employing the decomposition \eqref{Ue-dec} and taking $\boldsymbol\varphi^{\ep}$ as a test function in \eqref{H1-weak2} gives
\ba\label{limit-1}
& \int_{\Omega_\ep} \mu(\Theta) \left( \Grad \big(\vc{U}^0 - \vc{U}_{\ep,1}^0 - \vc{U}_{\ep,2}^0\big) + \Grad^t \big(\vc{U}^0 - \vc{U}_{\ep,1}^0 - \vc{U}_{\ep,2}^0\big) \right): \Grad \left(\boldsymbol\varphi - \boldsymbol\varphi^{\ep}_1 - \boldsymbol\varphi^{\ep}_2\right) \dx \\
& + \int_{\Omega_\ep} \mu(\Theta) \left( \Grad \vc{U}_\ep^r + \Grad^t \vc{U}_\ep^r \right): \Grad \boldsymbol\varphi^{\ep}\,\dx =  \int_{\Omega_\ep}   P_\ep \Div \left(\boldsymbol\varphi - \boldsymbol\varphi^{\ep}_1 - \boldsymbol\varphi^{\ep}_2\right)\dx + \langle \vc{f}_\ep, \boldsymbol\varphi^{\ep} \rangle_{W^{1,2}_0(\Omega_\ep)}.
\ea

\medskip

We first look at the right-hand side of \eqref{limit-1}.  By the convergence in \eqref{varphi-e-1-2-limit} and the assumption \eqref{P1-ass}, we have
$$
\lim_{\ep \to 0}  \langle \vc{f}, \boldsymbol\varphi^{\ep} \rangle_{W^{1,2}_0(\Omega_\ep)} = \langle \vc{f}, \boldsymbol\varphi \rangle_{W^{1,2}_0(\Omega)}
$$
and
$$
\limsup_{\ep \to 0} | \langle \vc{f}_\ep - \vc{f}, \boldsymbol\varphi^{\ep} \rangle_{W^{1,2}_0(\Omega_\ep)} |\leq C M \|\boldsymbol\varphi \|_{W^{1,2}_0(\Omega)}.
$$
Thus
\be\label{limit-fe}
\lim_{\ep \to 0} \langle \vc{f}_\ep, \boldsymbol\varphi^{\ep} \rangle_{W^{1,2}_0(\Omega_\ep)} = \langle \vc{f} , \boldsymbol\varphi \rangle_{W^{1,2}_0(\Omega)} + \langle \vc{r}, \boldsymbol\varphi \rangle_{W^{1,2}_0(\Omega)},
\ee
for some $\vc{r} \in W^{-1,2}_0(\Omega; R^3)$ satisfying
\[
\| \vc{r} \|_{W^{-1,2}_0 (\Omega; R^3)} \leq C M.
\]

\medskip

By the weak convergence of $P_\ep$ in \eqref{Ue-Pe-limit} and the strong convergence of $\boldsymbol\varphi^{\ep}_1$ in \eqref{varphi-e-1-2-limit}, we have
\ba\label{limit-Pe-1}
\int_{\Omega_\ep}   P_\ep \Div \boldsymbol\varphi \,\dx \to \int_{\Omega}   P \Div \boldsymbol\varphi \,\dx, \ \int_{\Omega_\ep}   P_\ep \Div \boldsymbol\varphi^{\ep}_1 \,\dx \to 0, \ \mbox{as $\ep \to 0$}.
\ea
By the definition of $\boldsymbol\varphi^{\ep}_2$ in \eqref{varphi-e-1-2}, and the divergence free property of $\vc{v}^i_{\ep}$, we have
$$
\Div \boldsymbol\varphi^{\ep}_2(x) =  \sum_{k=1}^{K(\ep)} \sum_{i=1}^3  \boldsymbol\varphi^i(x_k)  \vc{v}^i_{\ep} ( x - x_k ) \cdot \Grad \phi_{k,\ep}(x).
$$
Thus, by the property of $\vc{v}^i_{\ep}$ shown in Lemma \ref{model-pt} and the property of the cut-off function $\phi_{k,\ep}(x)$ in \eqref{chi-def}--\eqref{phi-k-e-def}, we have
\ba\label{limit-Pe-2}
\int_{\Omega_\ep}   P_\ep \Div \boldsymbol\varphi^{\ep}_2 \,\dx & = \int_{\Omega_\ep}   P_\ep \sum_{k=1}^{K(\ep)} \sum_{i=1}^3  \boldsymbol\varphi^i(x_k)  \vc{v}^i_{\ep} ( x - x_k ) \cdot \Grad \phi_{k,\ep}(x) \,\dx \\
& =  \sum_{k=1}^{K(\ep)} \sum_{i=1}^3 \int_{\{\frac{3}{4}\ep^3\leq |x-x_k| \leq \ep^3\}}   P_\ep   \boldsymbol\varphi^i(x_k)  \vc{v}^i_{\ep} ( x - x_k ) \cdot \Grad \phi_{k,\ep}(x) \,\dx \\
& \leq C \ep^{-3}  \int_{\bigcup_{k=1}^{K(\ep)} \{\frac{3}{4}\ep^3\leq |x-x_k| \leq \ep^3\}}  | P_\ep | \,\dx \\
& \leq C   \int_{\bigcup_{k=1}^{K(\ep)} \{\frac{3}{4}\ep^3\leq |x-x_k| \leq \ep^3\}}  | P_\ep |^2 \,\dx \\
& \to 0, \quad \mbox{as $\ep \to 0$}.
\ea

\medskip

We turn to consider the left-hand side of \eqref{limit-1}. By \eqref{Ur-limit3}, we have
\ba\label{limit-Ur}
 \int_{\Omega_\ep} \mu(\Theta) \left( \Grad \vc{U}_\ep^r + \Grad^t \vc{U}_\ep^r \right): \Grad \boldsymbol\varphi^{\ep}\,\dx \leq C h(\kappa),
\ea
which tends to $0$ when $\kappa \to 0.$

\medskip

By the  strong convergence of $\boldsymbol\varphi^{\ep}_1$ and weak convergence of $\boldsymbol\varphi^{\ep}_2$ in \eqref{varphi-e-1-2-limit}£¬ and similar convergence for $\vc{U}_{\ep,1}^0$ and $\vc{U}_{\ep,2}^0$, there holds
\ba\label{limit-U0-1}
& \int_{\Omega_\ep} \mu(\Theta) \left( \Grad \big(\vc{U}^0  - \vc{U}_{\ep,1}^0 - \vc{U}_{\ep,2}^0\big) + \Grad^t \big(\vc{U}^0 - \vc{U}_{\ep,1}^0 - \vc{U}_{\ep,2}^0\big) \right): \Grad \left(\boldsymbol\varphi - \boldsymbol\varphi^{\ep}_1 \right) \dx \\
 & \to \int_{\Omega} \mu(\Theta) \left( \Grad \vc{U}^0  +  \Grad^t \vc{U}^0  \right): \Grad \boldsymbol\varphi \,  \dx
\ea
and
\ba\label{limit-U0-2}
 \int_{\Omega_\ep} \mu(\Theta) \left( \Grad \big(\vc{U}^0 - \vc{U}_{\ep,1}^0 \big) + \Grad^t \big(\vc{U}^0 - \vc{U}_{\ep,1}^0 \big) \right): \Grad \boldsymbol\varphi^{\ep}_2 \, \dx \to 0
\ea
as $\ep \to 0$.

\medskip

It is left to study the limit
\ba\label{limit-U0-3}
\lim_{\ep \to 0} \int_{\Omega_\ep} \mu(\Theta) \left( \Grad \vc{U}_{\ep,2}^0  + \Grad^t  \vc{U}_{\ep,2}^0 \right): \Grad \boldsymbol\varphi^{\ep}_2 \, \dx.
\ea

By the definition in \eqref{varphi-e-1-2}, $\Grad \boldsymbol\psi^{\ep} $, where $ \boldsymbol\psi^{\ep} \in \{\vc{U}_{\ep,2}^0, \boldsymbol\varphi^{\ep}_2\}$, has two parts:
\be\label{psi-1}
\Grad \boldsymbol\psi^{\ep}  = \boldsymbol\psi^{\ep,1} +  \boldsymbol\psi^{\ep,2},
\ee
where
\be\label{psi-2}
\boldsymbol\psi^{\ep,1} : = \sum_{k=1}^{K(\ep)} \sum_{i=1}^3  \boldsymbol\psi^i(x_k)  \vc{v}^i_{\ep} ( x - x_k ) \Grad \phi_{k,\ep}(x) , \quad \boldsymbol\psi^{\ep,2} : = \sum_{k=1}^{K(\ep)} \sum_{i=1}^3  \boldsymbol\psi^i(x_k) \Grad \vc{v}^i_{\ep} ( x - x_k ) \phi_{k,\ep}(x).
\ee

By using \eqref{psi-1}--\eqref{psi-2}, we then can write \eqref{limit-U0-3} into four parts, and by a similarly argument as \eqref{limit-Pe-2}, any part involves $\boldsymbol\psi^{\ep,2}$ convergence to $0$ as $\ep \to 0$. Thus
\ba\label{limit-U0-4}
&\lim_{\ep \to 0} \int_{\Omega_\ep} \mu(\Theta) \left( \Grad \vc{U}_{\ep,2}^0  + \Grad^t  \vc{U}_{\ep,2}^0 \right): \Grad \boldsymbol\varphi^{\ep}_2 \, \dx\\
& = \lim_{\ep \to 0} \int_{\Omega_\ep} \mu(\Theta)  \sum_{k=1}^{K(\ep)} \sum_{i=1}^3  \vc{U}^{0,i}(x_k) \big(\Grad \vc{v}^i_{\ep} + \Grad^t \vc{v}^i_{\ep} \big) ( x - x_k ) \phi_{k,\ep}(x)  : \sum_{l=1}^{K(\ep)} \sum_{j=1}^3  \boldsymbol\varphi^j(x_l) \Grad \vc{v}^j_{\ep} ( x - x_l ) \phi_{l,\ep}(x)\,\dx \\
& = \lim_{\ep \to 0} \int_{\Omega_\ep} \mu(\Theta)  \sum_{k=1}^{K(\ep)} \sum_{i,j=1}^3  \vc{U}^{0,i} (x_k) \big(\Grad \vc{v}^i_{\ep} + \Grad^t \vc{v}^i_{\ep} \big) ( x - x_k ) \phi_{k,\ep}(x)  :   \boldsymbol\varphi^j(x_k) \Grad \vc{v}^j_{\ep} ( x - x_k ) \phi_{k,\ep}(x)\,\dx .
\ea

Again by a similarly argument as \eqref{limit-Pe-2}, and by the divergence free property of $\vc{v}^j_{\ep}$, we deduce from \eqref{limit-U0-4} that
\ba\label{limit-U0-5}
& \lim_{\ep \to 0} \int_{\Omega_\ep} \mu(\Theta)  \sum_{k=1}^{K(\ep)} \sum_{i,j=1}^3  \vc{U}^{0,i} (x_k) \Grad^t \vc{v}^i_{\ep} ( x - x_k ) \phi_{k,\ep}(x)  :   \boldsymbol\varphi^j(x_k) \Grad \vc{v}^j_{\ep} ( x - x_k ) \phi_{k,\ep}(x)\,\dx \\
& = \lim_{\ep \to 0} \int_{\Omega_\ep} \mu(\Theta)  \sum_{k=1}^{K(\ep)} \sum_{i,j=1}^3  \vc{U}^{0,i} (x_k)  \vc{v}^i_{\ep} ( x - x_k ) \phi_{k,\ep}(x)  \cdot   \boldsymbol\varphi^j(x_k) \Grad \vc{v}^j_{\ep} ( x - x_k ) \Grad^t \phi_{k,\ep}(x)\,\dx\\
& =0.
\ea

By \eqref{limit-U0-4}--\eqref{limit-U0-5}, by the definition of the cut-off function $\phi_{k,\ep}$ in \eqref{chi-def}--\eqref{phi-k-e-def}, there holds
\ba\label{limit-U0-6}
&\lim_{\ep \to 0} \int_{\Omega_\ep} \mu(\Theta) \left( \Grad \vc{U}_{\ep,2}^0  + \Grad^t  \vc{U}_{\ep,2}^0 \right): \Grad \boldsymbol\varphi^{\ep}_2 \, \dx\\
& = \lim_{\ep \to 0}  \mu(\Theta)  \sum_{k=1}^{K(\ep)} \sum_{i,j=1}^3  \vc{U}^{0,i} (x_k)  \boldsymbol\varphi^j(x_k) \int_{|x-x_k| \leq \ep^3} \phi_{k,\ep}^2(x)\Grad \vc{v}^i_{\ep}:\Grad \vc{v}^j_{\ep} ( x - x_k ) \,\dx .
\ea

Hence, by the assumption \eqref{i15}, together with Lemma \ref{model-pt}, we  obtain
\ba\label{limit-U0-f}
\lim_{\ep \to 0} \int_{\Omega_\ep} \mu(\Theta) \left( \Grad \vc{U}_{\ep,2}^0  + \Grad^t  \vc{U}_{\ep,2}^0 \right): \Grad \boldsymbol\varphi^{\ep}_2 \, \dx & = \sum_{i,j=1}^3  \int_{\Omega} \mu(\Theta)  \mathbb{D}_{i,j}(x) \vc{U}^{0,i} (x)  \boldsymbol\varphi^j(x)\,\dx \\
& =  \int_{\Omega} \mu(\Theta)  \mathbb{D} \vc{U}^{0} \cdot  \boldsymbol\varphi\,\dx.
\ea

\medskip

The final step is to pass $\kappa \to 0$. Recall the fact $\|\vc{U} - \vc{U}^0\|_{W^{1,2}_0(\Omega)} \leq \kappa$.
Then, by summarizing the limits in \eqref{limit-fe}, \eqref{limit-Pe-1}, \eqref{limit-Pe-2}, \eqref{limit-Ur}, \eqref{limit-U0-1}, \eqref{limit-U0-2} and \eqref{limit-U0-f} and by passing $\kappa\to 0$, we deduce that
\ba\label{limit-f}
& \int_{\Omega} \mu(\Theta) \left( \Grad \vc{U} + \Grad^t \vc{U}\right): \Grad \boldsymbol\varphi \, \dx + \int_{\Omega} \mu(\Theta)  \mathbb{D} \vc{U} \cdot  \boldsymbol\varphi\,\dx =  \int_{\Omega}   P \Div\boldsymbol\varphi \, \dx + \langle \vc{f}+ \vc{r}, \boldsymbol\varphi \rangle_{W^{1,2}_0(\Omega)}.\nonumber
\ea
The proof of Proposition \ref{P1} is completed.

\section{Asymptotic limit}
\label{a}

Our ultimate goal is to perform the asymptotic limit in the \emph{evolutionary} system (\ref{i1}--\ref{i3}).

\subsection{Compactness in time of the velocities}

We start by showing compactness in time of the family $\{ \vue \}_{\ep > 0}$ of the velocity fields. Let $\bfphi \in \DC(\Omega; R^3)$,
$\Div \bfphi = 0$. Set
\[
\vc{f} = - \Delta \bfphi + \mathbb{D} \cdot \bfphi \in L^\infty(\Omega; R^3).
\]
Let $\bfphi_\ep$ be the unique solution of the Stokes problem
\[
- \Delta \bfphi_\ep + \Grad P_\ep = \vc{f} \ ,\ \Div \bfphi = 0 \ \mbox{in}\ \Ome, \ \bfphi \in W^{1,2}_0 (\Ome; R^3).
\]
In accordance with Proposition \ref{P1},
\[
\bfphi_\ep \to \bfphi \ \mbox{weakly in}\ W^{1,2}_0 (\Omega; R^3); \ \mbox{whence}\ \bfphi_\ep \to \bfphi \ \mbox{in}\ L^2(\Omega; R^3).
\]
Now, we have
\[
\intO{ \vre \vue \cdot \bfphi } = \intO{ \vre \vue \cdot (\bfphi - \bfphi_\ep) } =
\intOe{ \vre \vue \cdot \bfphi_\ep },
\]
where, by virtue of the bounds (\ref{p1}), (\ref{p4}),
\begin{equation} \label{a1}
{\rm ess} \sup_{t \in (0,T)} \intO{ \vre \vue \cdot (\bfphi - \bfphi_\ep) } \to 0 \ \mbox{as}\ \ep \to 0.
\end{equation}
In addition, using $\psi (t) \bfphi_\ep$, $\psi \in \DC(0,T)$ as a test function in the variational formulation of the momentum balance (\ref{i10}),
we deduce that the family
\begin{equation} \label{a2}
\left\{ t \mapsto \intOe{ \vre \vue \cdot \bfphi_\ep } \right\}_{\ep > 0} \ \mbox{is precompact in}\ C([0,T]).
\end{equation}
Combining (\ref{a1}), (\ref{a2}),
we conclude that
\begin{equation} \label{a3}
\left[ t \mapsto \intO{ \vre \vue \cdot \bfphi } \right]
\to \left[ t \mapsto \intO{ \vr \vu \cdot \bfphi } \right] \ \mbox{in}\ L^\infty(0,T)
\ \mbox{for any}\ \bfphi \in \DC(\Omega; R^3), \ \Div \bfphi = 0.
\end{equation}
Using the density of smooth compactly supported functions in $W^{1,2}_{0, {\rm div}}(\Omega; R^3)$ - the Sobolev space $W^{1,2}_0(\Omega; R^3)$ of solenoidal
vector fields - we deduce from (\ref{a3}) that
\begin{equation} \label{a4}
\vre \vue \to \vr \vu \ \mbox{in}\ L^q([0,T]; W^{-1,2}_{{\rm div}}(\Omega; R^3)).
\end{equation}

Thus, finally, relation (\ref{a4}), together with (\ref{p5}), (\ref{p6}), imply that
\[
\int_0^T \intO{ \vr |\vue|^2 } \to \int_0^T \intO{ \vr |\vu|^2 },
\]
yielding
\begin{equation} \label{a5}
\vue \to \vu \ \mbox{in}\ L^2((0,T) \times \Omega; R^3).
\end{equation}

\subsection{Strong convergence of the temperature}

In view of (\ref{p4}), we may assume
\[
\Theta_\ep \to \Theta \ \mbox{weakly in}\ L^2(0,T; W^{1,2}_0(\Omega))
\ \mbox{and weakly-(*) in}\ L^\infty ((0,T) \times \Omega),
\]
passing to a subsequence as the case may be. Moreover, by virtue of (\ref{a5}),
the limit $\Theta$ solves (\ref{i3}) in the sense of (\ref{i11}), with $\kappa_\ep = \kappa_f$. Note that
\begin{equation} \label{a6}
\kappa_\ep \to \kappa_f \ \mbox{weakly-(*) in}\ L^\infty(\Omega) \ \mbox{and in}\ L^1(\Omega).
\end{equation}

As a matter of fact, $\Theta$ being a solution of the limit problem with \emph{constant} heat conductivity coefficient enjoys more regularity than
$\Theta_\ep$, specifically,
\begin{equation} \label{a7}
\Theta \in L^\infty (0,T; W^{2,2} \cap W^{1,2}_0(\Omega)).
\end{equation}

Finally, writing
\[
\begin{split}
\| \Theta_\ep - \Theta \|^2_{W^{1,2}_0(\Omega)} &\leq C \intO{ \kappa_\ep |\Grad \Theta_\ep - \Grad \Theta|^2 } \leq C
\intO{ \left( \kappa_{\ep}  \Grad \Theta_\ep - \kappa_f \Grad \Theta \right) \cdot \left(  \Grad \Theta_\ep - \Grad \Theta \right) }\\
& + C \intO{ (\kappa_f - \kappa_\ep) \Grad \Theta \cdot \left(  \Grad \Theta_\ep - \Grad \Theta \right) } \\
&= C \intO{ (\vue - \vu ) \cdot \Grad F (\Theta_\ep - \Theta )} + C  \intO{ (\kappa_f - \kappa_\ep) \Grad \Theta \cdot \left(  \Grad \Theta_\ep - \Grad \Theta \right) }
\end{split}
\]
we deduce from (\ref{a5}--\ref{a7}) that
\begin{equation} \label{a8}
\Theta_\ep \to \Theta \ \mbox{in}\ L^2(0,T; W^{1,2}_0(\Omega)).
\end{equation}
Relation (\ref{a8}), together with (\ref{p44}), yields the final conclusion
\begin{equation} \label{a9}
\Theta_\ep \to \Theta \ \mbox{in}\ L^q(0,T; C^\nu (\Ov{\Omega})) \ \mbox{for all}\ 1 \leq q < \infty
\ \mbox{and some}\ \nu > 0.
\end{equation}
Note that $\nu$ in (\ref{a9}) is strictly smaller than its companion in (\ref{p44}).

\subsection{Asymptotic limit in the momentum equation}

Our ultimate goal is to perform the asymptotic limit in the momentum equation (\ref{i10}). To this end, we use the time regularization by means of a convolution with a family of regularization kernels $\chi_\delta = \chi_\delta(t)$,
\[
\chi_\delta (t) = \frac{1}{\delta} \chi \left( \frac{t}{\delta} \right),\
\chi \in \DC (-1,1), \ \chi \geq 0, \ \chi(-z) = \chi(z), \ \chi'(z) \leq 0 \ \mbox{for}\ z \geq 0,\
\int_{-1}^1 \chi(z) \ {\rm d}z  = 1.
\]
As we are interested only in the behavior of $\vue$ on compact subsets of $(0,T) \times \Ome$, this step can be performed rigorously by
considering $\chi_\delta (\tau - t) \phi (x)$, $\phi \in \DC(\Omega; R^3)$, $\Div \phi = 0$ as a test function in the weak formulation (\ref{i10}).
Denoting $[v]_\delta = \chi_\delta * v$ we get
\begin{equation} \label{a10}
\intOe{ \left[ \mu (\Theta_\ep) \left( \Grad \vue + \Grad^t \vue \right) \right]_\delta : \Grad \phi } =
\intOe{ \left[ \left( \vre \vue \otimes \vue \right) \right]_\delta : \Grad \phi } - \intOe{ \partial_t [\vre \vue]_\delta \cdot \phi }.
\end{equation}
at any fixed $\tau \in (\delta, T - \delta)$.

Now, in view of (\ref{p5}), (\ref{p6}), and (\ref{a5}), it is easy to show that
\begin{equation} \label{a11}
[ \vre \vue \otimes \vue (\tau, \cdot) ]_{\delta} \to [\vr \vu \otimes \vu (\tau, \cdot) ]_{\delta} \ \mbox{in}\ L^2(\Omega; R^{3 \times 3}) \
\ \mbox{as}\ \ep \to 0,\  \mbox{for any}\
\tau \in (\delta, T - \delta),
\end{equation}
and, similarly,
\begin{equation} \label{a12}
\partial_t [\vre \vue (\tau, \cdot) ]_\delta \to \partial_t [\vr \vu (\tau, \cdot) ]_\delta
\ \mbox{in}\ L^2(\Omega; R^3) \ \mbox{as}\ \ep \to 0 \ \mbox{for any}\ \tau \in (\delta, T - \delta).
\end{equation}

Using (\ref{a11}), (\ref{a12}), we obtain the desired conclusion from (\ref{a9}) by application of Proposition \ref{P1} as soon as we show a suitable estimate
for the ``commutator''
\[
\left[ \mu (\Theta_\ep) \left( \Grad \vue + \Grad^t \vue \right) \right]_\delta -
\mu ( [\Theta]_\omega) \left( \Grad \left[ \vue \right]_\delta + \Grad^t \left[ \vue \right]_\delta \right),
\]
where $\delta > 0$, $\omega > 0$ are small parameters.

To begin, by virtue of (\ref{p4}), (\ref{a9}), observe that
\[
\begin{split}
\left[ \mu (\Theta_\ep) \left( \Grad \vue + \Grad^t \vue \right) \right]_\delta &-
\left[ \mu (\Theta) \left( \Grad \vue + \Grad^t \vue \right) \right]_\delta
\to 0 \ \mbox{in} \ L^2(\Omega; R^{3 \times 3}) \ \mbox{as}\ \ep \to 0 \\
&\mbox{uniformly for}\ \tau \in (\delta, T-\delta);
\end{split}
\]
whence
it is enough to control
\begin{equation} \label{a13}
\left[ \mu (\Theta) \left( \Grad \vue + \Grad^t \vue \right) \right]_\delta -
\mu ([\Theta]_\omega) \left( \Grad \left[ \vue \right]_\delta + \Grad^t \left[ \vue \right]_\delta \right).
\end{equation}

To handle (\ref{a13}), we write
\[
\begin{split}
 & \left[ \mu (\Theta) \left( \Grad \vue + \Grad^t \vue \right) \right]_\delta -
\mu ([\Theta]_\omega) \left( \Grad \left[ \vue \right]_\delta + \Grad^t \left[ \vue \right]_\delta \right) \\
&=\left[ (\mu (\Theta) - \mu([\Theta]_\omega) \left( \Grad \vue + \Grad^t \vue \right) \right]_\delta + \left[ \mu ([\Theta]_\omega) \left( \Grad \vue + \Grad^t \vue \right) \right]_\delta\\
 &\quad - \mu ([\Theta]_\omega) \left( \Grad \left[ \vue \right]_\delta + \Grad^t \left[ \vue \right]_\delta \right).
\end{split}
\]
Now,
\[
\begin{split}
&\left\|
\left[ (\mu (\Theta) - \mu([\Theta]_\omega) \left( \Grad \vue + \Grad^t \vue \right) \right]_\delta \right\|_{L^2(\Omega; R^{3 \times 3})} \\ &\leq
\left\|
(\mu (\Theta) - \mu([\Theta]_\omega) \left( \Grad \vue + \Grad^t \vue \right)  \right\|_{L^2(\Omega; R^{3 \times 3})}\\
& \leq
\left\|
\mu (\Theta) - \mu([\Theta]_\omega) \right\|_{L^\infty(\Omega)} \left\| \Grad \vue + \Grad^t \vue   \right\|_{L^2(\Omega; R^{3 \times 3})}.
\end{split}
\]

Thus, in view of (\ref{p4}), (\ref{a9}),
\begin{equation} \label{a14}
\left[ (\mu (\Theta) - \mu([\Theta]_\omega) \left( \Grad \vue + \Grad^t \vue \right) \right]_\delta
\to 0 \ \mbox{in}\ L^q(0,T; L^2(\Omega; R^{3 \times 3})) \to 0 \ \mbox{as}\ \omega \to 0
\end{equation}
for any $1 \leq q < 2$,
uniformly in $\ep$ and $\delta$.

On the other hand, if $\omega > 0$ is fixed, the function $[\Theta]_\omega$ is continuously differentiable with respect to the spatial variable.
Thus we deduce that
\begin{equation} \label{a15}
\begin{split}
&\left[ \mu ([\Theta]_\omega) \left( \Grad \vue + \Grad^t \vue \right) \right]_\delta -
\mu ([\Theta]_\omega) \left( \Grad \left[ \vue \right]_\delta + \Grad^t \left[ \vue \right]_\delta \right)\\
&= \left[ \mu ([\Theta]_\omega) \left( \Grad \vue + \Grad^t \vue \right) \right]_\delta -
\mu ([\Theta]_\omega) \left[ \Grad  \vue  + \Grad^t \vue \right]_\delta \\
&\to 0
\ \mbox{in}\ L^2((0,T) \times \Omega; R^{3 \times 3}) \ \mbox{as}\ \delta \to 0
\end{split}
\end{equation}
uniformly in $\ep$ for any fixed $\omega > 0$.

Summing up relations (\ref{a14}), (\ref{a15}), we may rewrite (\ref{a10}) in the form
\[
\begin{split}
\intOe{ \mu ([\Theta]_\omega) \left( \Grad [\vue]_\delta + \Grad^t [\vue]_\delta \right)  : \Grad \phi } &=
\intOe{ \left[ \left( \vre \vue \otimes \vue \right) \right]_\delta : \Grad \phi } - \intOe{ \partial_t [\vre \vue]_\delta \cdot \phi }\\
&+ \intOe{ \left( \mathbb{R}^1_{\omega, \delta, \ep} + \mathbb{R}^2_{\omega, \delta,\ep } \right) : \Grad \phi }
\end{split}
\]
at any fixed $\tau \in (\delta, T - \delta)$, where
\[
\begin{split}
\mathbb{R}^1_{\omega, \delta, \ep} &\to 0 \ \mbox{in}\ L^q(0,T; L^2(\Omega; R^{3 \times 3}))
\ \mbox{as}\ \omega \to 0,\ 1 \leq q < 2, \ \mbox{uniformly in}\ \ep, \delta,\\
\mathbb{R}^2_{\omega, \delta, \ep} &\to 0 \ \mbox{in}\ L^2(0,T; L^2(\Omega; R^{3 \times 3}))
\ \mbox{as}\ \delta \to 0 \ \mbox{for any fixed}\ \omega > 0 \ \mbox{and uniformly in}\ \ep.
\end{split}
\]

Thus performing successively the limits $\ep \to 0$, $\delta \to 0$, and, finally, $\omega \to 0$, we deduce the desired conclusion. Theorem \ref{T1} has been proved.

\def\cprime{$'$} \def\ocirc#1{\ifmmode\setbox0=\hbox{$#1$}\dimen0=\ht0
  \advance\dimen0 by1pt\rlap{\hbox to\wd0{\hss\raise\dimen0
  \hbox{\hskip.2em$\scriptscriptstyle\circ$}\hss}}#1\else {\accent"17 #1}\fi}


\begin{thebibliography}{1}

\bibitem{Allai4}
G.~Allaire.
\newblock Homogenization of the {N}avier-{S}tokes equations in open sets
  perforated with tiny holes. {I}. {A}bstract framework, a volume distribution
  of holes.
\newblock {\em Arch. Rational Mech. Anal.}, 113(3), 209--259, 1990.

\bibitem{Allai3}
G.~Allaire.
\newblock Homogenization of the {N}avier-{S}tokes equations in open sets
  perforated with tiny holes. {II}. {N}oncritical sizes of the holes for a
  volume distribution and a surface distribution of holes.
\newblock {\em Arch. Rational Mech. Anal.}, 113(3), 261--298, 1990.

\bibitem{CHAN}
S.~Chandrasekhar.
\newblock {\em Hydrodynamic and hydrodynamic stability}.
\newblock Clarendon Press, Oxford, 1961.

\bibitem{CioMur2}
D.~Cioranescu and F.~Murat.
\newblock Un terme \'etrange venu d'ailleurs.
\newblock In {\em Nonlinear partial differential equations and their
  applications. {C}oll\`ege de {F}rance {S}eminar, {V}ol. {II} ({P}aris,
  1979/1980)}, volume~60 of {\em Res. Notes in Math.}, pages 98--138, 389--390.
  Pitman, Boston, Mass., 1982.

\bibitem{CioMur1}
Do{\"{\i}}na Cioranescu and Fran{\c{c}}ois Murat.
\newblock Un terme \'etrange venu d'ailleurs. {II}.
\newblock In {\em Nonlinear partial differential equations and their
  applications. {C}oll\`ege de {F}rance {S}eminar, {V}ol. {III} ({P}aris,
  1980/1981)}, volume~70 of {\em Res. Notes in Math.}, pages 154--178,
  425--426. Pitman, Boston, Mass., 1982.


\bibitem{DesGolRic}
L.~Desvillettes, F.~Golse, and V.~Ricci.
\newblock The meanfield limit for solid particles in a {N}avier--{S}tokes flow.
\newblock {\em J. Stat. Physics}, {\bf 131}, 941--967, 2008.

\bibitem{DFL17} L. Diening, E. Feireisl, Y. Lu. The inverse of the divergence operator on perforated domains with applications to homogenization problems for the compressible Navier--Stokes system. {\em ESAIM: Control Optim. Calc. Var.,}  {\bf 23} (3) (2017), 851-868.


\bibitem{DL}
R.J. DiPerna and P.-L. Lions.
\newblock Ordinary differential equations, transport theory and {S}obolev
  spaces.
\newblock {\em Invent. Math.}, {\bf 98}:511--547, 1989.

\bibitem{FNN16} E. Feireisl, Y. Namlyeyeva,  \v{S}\'arka Ne\v{c}asov\'a. Homogenization of the evolutionary Navier--Stokes system. Manuscripta Math. 149, 251--27 (2016)


\bibitem{FENO6}
E.~Feireisl and A.~Novotn{\' y}.
\newblock {\em Singular limits in thermodynamics of viscous fluids}.
\newblock Birkh{\" a}user-Verlag, Basel, 2009.


\bibitem{LADUR}
O.~A. Ladyzhenskaya and N.~N. Uralceva.
\newblock {\em Linear and quasilinear elliptic equations}.
\newblock Academic Press, New York and London, 1968.

\bibitem{LIGN}
F.~Ligni\`eres.
\newblock The small {P}\'eclet number approximation in stellar radiative zones.
\newblock {\em Astronomy and astrophysics}, pages 1--10, 1999.

\bibitem{Lions-Incom} P.-L. Lions. \newblock {\em Mathematical topics in fluid dynamics, Vol.1, Incompressible  models}. \newblock Oxford Science Publication, Oxford, 1998.

\bibitem{MarKhr}
V.~A. Marchenko and E.~Ya. Khruslov.
\newblock {\em Homogenization of partial differential equations}, volume~46 of
  {\em Progress in Mathematical Physics}.
\newblock Birkh\"auser Boston, Inc., Boston, MA, 2006.
\newblock Translated from the 2005 Russian original by M. Goncharenko and D.
  Shepelsky.

\end{thebibliography}

\end{document}